\documentclass{article}
\catcode`\@=11
\@addtoreset{equation}{section}

\catcode`\@=12\newtheorem{Theorem}{Theorem}[section]
\newtheorem{Proposition}{Proposition}[section]
\newtheorem{Lemma}{Lemma}[section]

\usepackage{graphicx}

\font\fil=msbm10 
\def\R{\mbox{\fil R}}\def\N{\mbox{\fil N}}
\def\Rn{\R^n}
\def\A{{\cal A}}
\def\qed{\vbox{\hrule\hbox to7.8pt{\vrule height7pt
    \hss\vrule height7pt}\hrule}}

\newcommand{\cl}{\hbox{cl}}

 \newcommand\DomPhi{\hbox{D}_\Phi}
\title{Global Inversion of Functions: \\
an Introduction
\thanks{This paper appeared on \emph{NoDEA} 1  (1994) 229--248.
 Authors' email addresses in 2014: giuseppe.demarco@unipd.it,
 gianluca.gorni@uniud.it, gaetano.zampieri@univr.it}}
\author{Giuseppe De Marco\\
Universit\`a di Padova,  Dip. di
 Matematica Pura e
Appl.\\
 via Belzoni 7, 35131 Padova, Italy\\
\and
Gianluca Gorni\\
 Universit\`a di Udine,
Dip. di Matematica e Informatica\\
 via Zanon 6, 33100 Udine, Italy\\
 email: Gorni@udmi5400.cineca.it
\and
Gaetano Zampieri\\
Universit\`a di Padova,  Dip. di
 Matematica Pura e
Appl.\\
 via Belzoni 7, 35131 Padova, Italy\\
email: Gaetano@pdmat1.unipd.it}
\date {{\sl Dedicated to Roberto Conti on the occasion of
 his 70th birthday}}
\begin{document}
\textwidth=125mm
\textheight=185mm
\parindent=8mm
\frenchspacing
\maketitle 
\begin{abstract}This is an exposition of some basic ideas in 
the realm of Global Inverse
Function theorems. We   address ourselves mainly to
readers who are interested 
in the applications to Differential Equations.
But we do not deal with those applications and
we  give a `self-contained' elementary exposition.

The first part is devoted to
the celebrated Hadamard-Caccioppoli theorem on
proper local homeomorphisms
treated in the framework of the Hausdorff spaces.
In the proof, the concept of
`$\omega$-limit set' is used in a  crucial 
way and this is perhaps the novelty
of our approach. 

In the second part we deal with open sets 
in Banach spaces.
The concept of `attraction basin'  here is
the main tool of our exposition
which  also shows a few recent results, here extended 
from finite dimensional
to general Banach spaces, 
together with the classical theorem 
of Hadamard-Levy which assumes that the
operator norm of the
inverse of the derivative does not grow
too fast (roughly at most linearly).
\end{abstract}

\section*{Introduction}
A fundamental problem in Analysis is the existence and/or
uniqueness of the solutions to the equation $y=f(x)$
in the unknown $x$. The function $f:X\to Y$ relates
two spaces $X,Y$ with some structure, otherwise we are
impotent. From the other side, the concrete  case
where $X,Y$ are subsets of the n-space $\Rn$ is often too
restrictive,
and actually many applications arise
in more general  spaces. We especially think
about injectivity and surjectivity problems in Differential
Equations which are not discussed in this paper but constitute
one of the reasons of our discussion. 

The books  Prodi
and Ambrosetti [31], and Chow and Hale [9], give the proof 
of global inversion theorems in general spaces and show
applications to differential equations.  
Let us also refer to Invernizzi and Zanolin [21], Brown and Lin
[6], and  Radulescu and
Radulescu [33]  among the   papers
which could be mentioned for results in differential equations
obtained by means of the inversion of functions in infinite
dimensional Banach spaces.
Finite dimensional problems are also important.
The research field of the Jacobian conjectures deals with
deep questions
of invertibility linked to global stability problems, see Olech
[27], Meisters [23], Meisters and Olech [25], [26],
and the references
contained therein. The inversion of functions, of course,
also plays a role  in the applied sciences, e.g.  Economics
and Network Theory.

More references are listed at the end of the paper with no
claim to completeness. The present paper is not a survey on
the rich literature on these topics.

Section 1 below is devoted to the following theorem which
we call after Hadamard and Caccioppoli since Hadamard
was probably the first  to have the idea in finite
dimension, and Caccioppoli was perhaps the most important
author in the process of clarification and
generalization to abstract spaces (but 
other mathematicians also gave a contribution).
\begin{Theorem}[Hadamard-Caccioppoli].  Let
$f:X\to Y$ be a local homeomorphism with $X,Y$ path
connected Hausdorff spaces and $Y$ simply connected. Then $f$
is a homeomorphism onto $Y$ if and only if it is a proper
function, namely if and only if the inverse image
$f^{\leftarrow}(K)$ of any compact set $K\subset Y$ is
compact.\end{Theorem}
The proof below  uses, in a crucial way, the concept of
$\omega$-limit set.
This is perhaps the main novelty of our approach. 

The statements of the Theorem in the books of Prodi and
Ambrosetti [31], and Chow and Hale [9]
(whose treatment of this topic is based on [31]), seem different from
Theorem 1 at a first glance since they mention possible  singular
points of $f$; however those statements actually follow at once
from the one above. Incidentally, those books state
the theorem in metrizable spaces. 
We believe that the more general framework of Hausdorff
spaces does not cost more than usual
presentations in metrizable 
spaces even if these are, of course, the
relevant case for applications. And  generality usually favours  
 understanding  the essence of a subject. The framework 
of Theorem 1 is somehow essential, in particular it is false in
non-Hausdorff topological spaces  as a simple counterexample
will show.

Finally we show an application of the theorem to Algebra, due
to Gordon. Namely we show, following [14],
that there cannot be a product in $\Rn$
for $n\geq3$ (see Proposition~1.3 below for a precise
formulation). This is related to the fact that $\Rn
\setminus\{0\}$ is simply-connected if and only if $n\geq3$.
We  quote this application to convince the reader of the depth
of the Hadamard-Caccioppoli theorem in a concise way.

In Section 2 we deal with local homeomorphisms $f:D\to Y$ from
an open connected set of a Banach space $X$ to a Banach space
$Y$.  In order to briefly mention the ideas discussed there, let
us here refer to the particular case of a local diffeomorphism
$f$. Then the celebrated
 Wa\.zewski equation with parameter $v\in Y$, 
  
  \begin{equation}
\dot x=f'(x)^{-1}\,v
\end{equation}  
  is often used in the literature to deal with invertibility
  problems.
  Wa\.zewski  introduced (0.1) in  [40], for $X=Y=\Rn$,
  to give an estimate for a ball, around a given point
  $x_0\in D$, where the inverse
  function can be defined. Instead of (0.1) we consider
  \begin{equation}
\dot x=F(x)\,,\qquad F: D\to X\,,\  x\mapsto
  -f'(x)^{-1}\left(f(x)-f(x_0)\right)
\end{equation}
  whose trajectories are
  also trajectories of the family of equations (0.1) (as $v\in
  Y$) but with different
  parametrization (incidentally, remark that the family (0.1)
  has many more trajectories). 
  
  The point $x_0$ is an asymptotically stable equilibrium
  for (0.2) and its attraction basin $\A$ will be proved to 
  coincide with the maximal open
  subset of $D$, containing $x_0$,
  such that $f|\A$ is injective and, at the same
  time, the image $f(\A)$ is star-shaped with respect to
  $y_0:=f(x_0)$. Using these ideas we show some criteria for
  the injectivity of $f$. Moreover, we shall see that
  the solutions to the equation (0.2) are all defined on 
  the whole $\R$ if and only if $f$ is a global homeomorphism
  onto $Y$. In particular, this fact leads to the following:
 \begin{Theorem}[Hadamard--Lev\'y]Let
  $f:X\to Y$ be a local diffeomorphism with $X,Y$ 
  Banach spaces. Then $f$
  is a diffeomorphism onto $Y$ if there exists a continuous
  (weakly) increasing 
  map $\beta:\R_+ \to \R_+\setminus \{0\}$ such that
 
\begin{equation}
\int_0^{+\infty}\;\frac1{\beta(s)}ds=
+\infty\,,\qquad
  \|f'(x)^{-1}\|\leq\beta(\|x\|).
\end{equation}
  In particular this holds if, for some $a,b\in\R_+$, we have
  \begin{equation}\|f'(x)^{-1}\|\leq
a+b\|x\|.\end{equation}\end{Theorem}
 This theorem was discovered by Hadamard in $\R^n$.
  Then it was generalized by Levy to infinite dimension 
  under condition (0.4) with\  $b=0$. Meyer dealt with the
  full condition (0.4), and finally Plastock gave a proof for
  the general statement. In the literature it is often named
  after Hadamard only. 
  
  Finally, we deal with the  injectivity of $f$ 
  (together with the star-shape of the image) by means of
  global Lyapunov functions. We extend to general Banach spaces 
  some  results
  previously obtained in [17] by two of the authors for $\R^n$.
  
  Our approach to the invertibility of functions, by means
  of attraction basins for (0.2), is one of the
  ingredients used in
  [26] by Meisters and Olech
  to prove one of the results in that paper,
  namely the global asymptotic stability 
  for certain polynomial vector fields.
  We hope that it can lead
  to further consequences, in particular for the
  Differential Equations.

\section{The Hadamard-Caccioppoli Theorem}
  In this Section $X,Y,Z$ will always be topological
Haus\-dorff spaces.

\medskip
 {\bf Local homeomorphism}. As is well known the
function $f:X\to Y$ is called  a local homeomorphism at
$x_0\in X$ if there exist open neighbourhoods $U,V$ of $x_0$
and $y_0:=f(x_0)$ respectively, such that $f(U)=V$ and the
restriction $f|U:U \to V$ is a homeomorphism.
Then $g:=(f|U)^{-1}:V\to U$ is called a local inverse of $f$
at $y_0$. Moreover we say that $f:X\to Y$ is a local
homeomorphism if it is a local homeomorphism at any $x_0\in
X$. Such a mapping is clearly continuous and open, namely 
 inverse-images and images of open sets are open sets.

\medskip

{\bf Lifting}. Let $f:X\to Y$ be a local
homeomorphism and let $p :Z\to Y$ be a continuous
function. A continuous function $\tilde{p} : Z\to X$ is
called a lifting of $p$ by $f$ whenever $f\circ
\tilde{p} = p$, that is if the following diagram
commutes:

\medskip

\[ 
\begin{array}{rcr}
\null&\null&X\\
\noalign{\medskip}
\null&{}\raise5pt\rlap{$\scriptstyle\tilde p$}\nearrow 
&{}\downarrow\rlap{$\scriptstyle f$}\\
\noalign{\medskip}
Z&\stackrel{p}{\longrightarrow}&Y
\end{array}
\]  

\medskip

  \begin{Lemma}(Uniqueness). 
  Let $f:X\to Y$ be
a local homeomorphism between Hausdorff spaces and let
$p:Z\to Y$ be continuous with $Z$ connected. If $\ \tilde
p_1,\tilde{p}_2:Z\to X$ are both liftings of
$p$ then either
$\tilde{p}_1=\tilde{p}_2$ or
$\tilde{p}_1(z)\ne\tilde{p}_2(z)$ for every $z\in
Z$. \end{Lemma}
{\it Proof.} Let
$C:=\{z\in Z:\tilde{p}_1(z)=\tilde{p}_2(z)\}$. Let
us see that $C$ is open in $Z$. If $C=\emptyset$ then it is
open; otherwise take $z_0\in C$ and let $x_0:=
\tilde{p}_1(z_0)=\tilde{p}_2(z_0)$. Moreover let
$U,V$ and $g:V\to U$ be as in the definition of local
homeomorphism above. The set
$W:=\tilde{p}_1^{\leftarrow}(U)\cap
\tilde{p}_2^{\leftarrow}(U)$ is an open neighbourhood of $z_0$
and we have $\tilde{p}_1|W=\tilde{p}_2|W=g\circ p|W$. Thus
$W\subseteq C$ and $C$ is open.

Now, $Z\setminus C$ is open 
by an easy standard argument (which uses  that $X$ is
Hausdorff), so we are done since $Z$ is connected. 

 \hfill\qed 
 
\bigskip\goodbreak
 {\bf Path-lifting property}. We say that the
local homeomorphism $f:X\to Y$ lifts the paths if, for every
continuous function $\alpha:[0,1]\to Y$, with $\alpha(0)\in
f(X)$ (called a path in $Y$ with origin in $f(X)$), and
for every $x_0\in f^{\leftarrow}(0)$, there exists a lifting
$\tilde\alpha:[0,1]\to X$ of $\alpha$ with
$\tilde{\alpha}(0)=x_0$. By Lemma 1.1, if $f$ lifts the
paths then it does it with uniqueness, that is the
$\tilde{\alpha}$ above is unique.

\medskip
{\bf Homotopy-lifting property}. A continuous map $H:Z\times
[0,1]\to Y$ is called a homotopy with base $H_0:Z\to Y$,
$z\mapsto H(z,0)$. We say that $f:X\to Y$ lifts the homotopies
if, for any such $H$, and any continuous 
map ${\tilde H}_0:Z\to X$ such that $f\circ {\tilde H}_0= H_0$
(${\tilde H}_0$ is a lifting of  the base of the homotopy),
there exists a continuous lifting $\tilde H$ with base ${\tilde
H}_0$, that is $f\circ \tilde H=H$ and ${\tilde H}(z,0)={\tilde
H}_0(z)$ for all $z\in Z$.

The path-lifting property is clearly a particular case of the
homotopy-lifting property, with $Z$ a one-point space. It is
then remarkable the following

 \begin{Lemma}
(Path-lifting $\Longrightarrow$ Homotopy-lifting). 
 If the local homeomorphism between Hausdorff spaces
 $f:X\to Y$ lifts the paths, then it lifts the
homotopies.\end{Lemma}

{\it Proof}.   With the notations as in the above definitions, 
let $t\mapsto \tilde H(z,t)$ be the unique lifting of the path
$t\mapsto H(z,t)$, with origin ${\tilde H}_0(z)$, for any $z\in
Z$. Clearly $f\circ\tilde H=H$, and $\tilde H(z,0)={\tilde
H}_0(z)$. So starting from $H$ and ${\tilde H}_0$ as above, we
have  defined  $\tilde H$, all we are left to prove is its
continuity   on $Z\times[0,1]$. Take $z_0\in Z$,
and let $D$ be the subset of $[0,1]$ consisting of all $t\in
[0,1]$ such that $\tilde H$ is not continuous at $(z_0,t)$. We
argue by contradiction: assuming $D$ non empty, $D$ has an
infimum $a\ge 0$; since $t\mapsto\tilde H(z_0,t)$ is
continuous, given any neighborhood $U$ of $\tilde H(z_0,t_0)$
in $X$ there exists an interval $J_1$, an open neighborhood of
$a$ in $[0,1]$, such that $\tilde H(z_0,t)\in U$ for every
$t\in J_1$. By restricting $U$ if necessary we can assume  $U$ 
open, and that $f$ induces a homeomorphism $f|U:U\to V$ onto
a neighborhood $V$ of $H(z_0,a)$. By continuity of $H$  there
exists a neighborhood $W_1$ of $z_0$ in $Z$, and another
interval  $J_2$, open neighborhood of $a$ in $[0,1]$, such
that $H(W_1\times J_2)\subseteq V$. Let $J=J_1\cap J_2$, and
pick   $b\in J$, with $b<a$ if $a>0$; if $a=0$
let $b=0$;  in both cases $z\mapsto\tilde H(z,b)$ is
continuous at $z_0$ (as a function from $Z$ to $X$), and since
$\tilde H(z_0,b)\in U$, with $U$ open, there exists a
neighborhood $W_2$ of $z_0$ in $Z$ such that $\tilde
H(W_2\times\{b\})\subseteq U$; put $W=W_1\cap W_2$. We
claim that $$\tilde H|W\times J=(f|U)^{-1}\circ H|W\times
J;$$ in fact these functions coincide on $W\times\{b\}$;
but then, for every $z\in W$ the functions defined on $J$ by
$t\mapsto\tilde H(z,t)$, $t\mapsto (f|U)^{-1}\circ H(z,t)$ are
liftings of $t\mapsto H(z,t)$ which coincide on $b\in J$,
an hence coincide on all of $J$. The equality just proved shows
that $\tilde H$ is continuous at $(z_0,t)$, for every $t\in
J,\quad t\ge a$,  contradicting the minimality of $a$.

\hfill\qed

\bigskip

\begin{Lemma}  (Simply connected
codomain).  Let $f:X\to Y$ be
a local homeomorphism between Hausdorff spaces which lifts the
paths. If $X,Y$ are path connected and $Y$ is simply
connected, then $f$ is a homeomorphism.\end{Lemma}

{\it Proof}. First of all let us see the surjectivity. Let
$y_0\in f(X)$, $x_0\in f^{\leftarrow}(y_0)$, and let $\alpha:
[0,1]\to Y$ be a path with $\alpha (0)=y_0$ and $\alpha
(1)=y$. There exists a (unique) lifting $\tilde\alpha$ of
$\alpha$ with $\tilde \alpha(0)=x_0$. The formula
$f\circ\tilde\alpha=\alpha$ gives $f(\tilde\alpha(1))=y$.

Now, let us see the injectivity of $f$. Let $x_0,x_1\in X$
satisfy
$f(x_0)=f(x_1)=:y_0$. Since $X$ is path connected we can
consider a path $\sigma:[0,1]\to X$ joining $x_0,x_1$, that is
with $\sigma(0)=x_0$ and $\sigma(1)=x_1$. The formula
$\alpha:=f\circ\sigma$ defines a circuit in $Y$ (i.e. a closed
path) with $\alpha(0)=\alpha(1)=y_0$. Since $Y$ is simply
connected there exists a homotopy with fixed end-points $h$
between $\alpha$ and the  constant path $[0,1]\to Y, t\mapsto
y_0$, namely a continuous function $h:[0,1]^2 \to Y$ such that
$h(t,0)=\alpha(t)$, $h(t,1)=y_0$ for all $t\in [0,1]$, and
$h(0,s)=y_0=h(1,s)$,  for all
$s\in [0,1]\,$  (see
the figure below). 

Since $f$ lifts  paths, then, by Lemma~1.2, there
exists a unique $\tilde h:[0,1]^2\to X$ which lifts $h$ and
which satisfies $\tilde h(t,0)=\sigma(t)$, for all $t\in[0,1]$.

In the rest of the proof we use the following important fact:
a constant path is lifted to a constant path (which works
being continuous and which is the unique lifting by
Lemma~1.1). Thus $\tilde h(0,s)=\sigma(0)=x_0$, $\tilde h(1,s)
=\sigma(1)=x_1$, for all $s\in[0,1]$; and since 
$t\mapsto\tilde h(t,1)$ is also constant, we have $x_0=\tilde h(0,1)
=\tilde h(1,1)=x_1$.

\hfill\qed

\bigskip\bigskip
\centerline{\includegraphics[scale=.5]{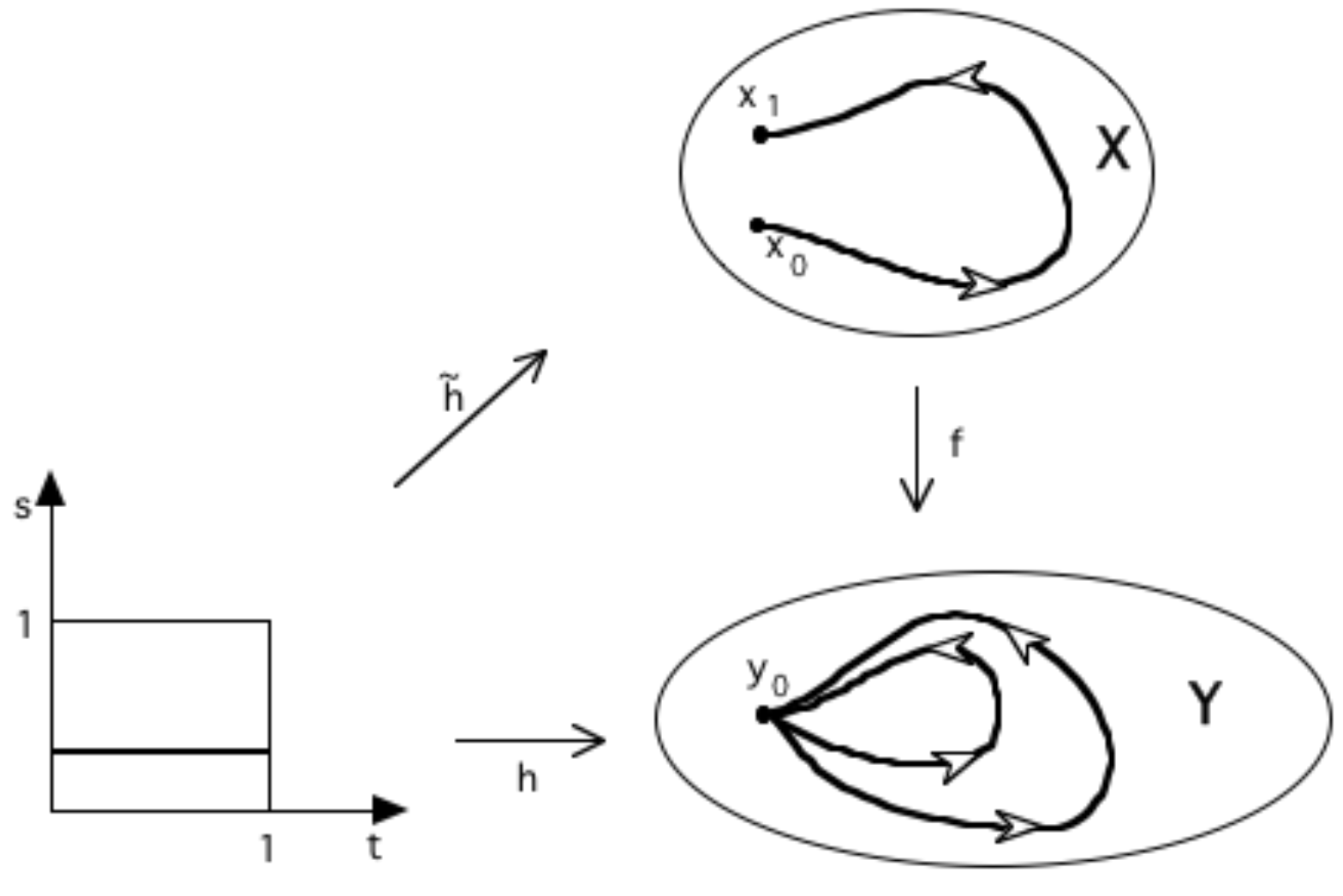}}
 
\bigskip\bigskip

\goodbreak

{\bf Maximal path-lifting}. Let $f:X\to Y$ be a local
homeomorphism,  let $\alpha:[0,1]\to Y$ be a path with
$\alpha(0)\in f(X)$, and let $x_0\in
f^{\leftarrow}(\alpha(0))$.
We  define the maximal lifting $\phi:J\to X$ of $\alpha$
with $\phi(0)=x_0$ in the following way. There certainly
exists a continuous map $\phi_I:I\to X$, with $I=[0,b[\subset
[0,1]$, such that $\phi_I(0)=x_0$ and $f\circ \phi_I=\alpha
|I$.  By the uniqueness Lemma 1.1, the formula $\phi
|I=\phi_I$ defines the mapping $\phi :J\to X$ on the union
$J$ of all the intervals $I$.

\medskip 
{\bf $\omega$-limit set}. Let $\phi:[0,b[\to X, \;
0<b\le +\infty$ be a continuous function. Then the following
formula, where `$\cl$' denotes the closure in $X$, defines the
$\omega$-limit set of $\phi$:
\[\omega_{\phi}:=\bigcap_{t\in [0,b[}\;
\cl\,\phi([t,b[\,)\,.\]

Equivalently, $x\in\omega_{\phi}$ if and only if $x$ is a
cluster point of a sequence $(\phi(t_n))$, for some sequence
$t_n\in[0,b[$ which converges to $b$; in the particular case
of $X$  metrizable, $x\in  \omega_{\phi}$ if and only
if there exists a sequence $(t_n)$ with $t_n\in [0,b[$ such
that  $t_n\to b$ and $\phi (t_n)\to x$ as $n\to \infty$.

If $\phi$ were a solution of an autonomous
differential equation $\dot x=F(x)$, then the terminology
`$\omega$-limit set' would be usual. This concept has 
paramount importance since one of the main goal of Dynamics is
precisely to say what is the destiny of the motions
(incidentally, recall that $\omega$ is the last letter of the
Greek alphabet). 
 
\goodbreak\bigskip

 \begin{Lemma}($\omega$-limit set
of a maximal path lifting).   Let $f:X\to Y$ be a local
homeomorphism between Hausdorff spaces, 
 and
let  $\phi:J\to X$ be the maximal lifting of $\alpha:[0,1]\to
Y$ with $\phi(0)=x_0\in f^{\leftarrow}(\alpha(0))$. 
If $J\ne [0,1]$ then it is open to the right, i.e. $J=[0,b[$ 
with  $b\in ]0,1]$, and the $\omega$-limit set of $\phi$ is
empty: $\omega_{\phi}=\emptyset$.\end{Lemma}

\goodbreak\bigskip

{\it Proof}. We argue by contradiction by assuming that
$J=[0, a]$ with $0<a<1$. We  consider a local inverse
of $f$ at $f(\phi(a))$ and we easily extend $\phi$ to a lifting
defined on a larger domain, this contradicts the maximality
of $\phi$. So  $\phi:[0,b[\to X$ for a suitable $b\in ]0,1]$.
 
 Now, let us contradict $\omega_{\phi}=\emptyset$ and let
 $x_0\in \omega_{\phi}$. Then $f(x_0)=\alpha (b)$ since 
 by  continuity
 $f(\cl\, \phi ([t,b[\,))\subseteq \cl\, f(\phi
([t,b[\,)$ and
 \[ 
 \bigcap_{t\in [0,b[} \cl\,f(\phi ([t,b[\,)=
 \bigcap_{t\in [0,b[} \alpha([t,b])=\{\alpha(b)\}\]
  (in metric spaces we could just argue with
 sequences).
 
Consider open neighbourhoods $U,V$, of $x_0$ and $f(x_0)$
respectively, such that $f|U:U\to V$ be a homeomorphism, and
let $g$ be the inverse function. We can consider $a\in [0,b[$
such that $\alpha ([a,b])\subset V$, and such that $\phi(a)\in
U$. Moreover, we can define $\psi :[0,b]\to X$ lifting of
$\alpha |[0,b]$ by $\psi|[0,a]=\phi |[0,a]$ and by
$\psi |]a,b]=g\circ \alpha |]a,b]$. This contradicts the
maximality of $\phi$.

\hfill \qed

\bigskip

Now we are ready to prove Theorem~0.1 of the Introduction.
\goodbreak\bigskip

{\bf Proof of the Hadamard-Caccioppoli Theorem}. Let $f$ be
proper (in the other sense the theorem is trivial). We are going
to prove that $f$ lifts the paths.  This gives the theorem  by
means of Lemma~1.3.

We argue by contradiction by assuming the existence of a path
$\alpha:[0,1]\to Y$ and a point $x_0\in
f^{\leftarrow}(\alpha(0))$ such that the maximal lifting $\phi$
of $\alpha$, with $\phi(0)=x_0$, is defined on $[0,b[$,
with $b\le 1$ (but not on $[0,1]$). Then Lemma~1.4 says
that $\omega_{\phi}=\emptyset$.

But $\phi([0,b[\,)\subset
f^{\leftarrow}(\alpha([0,1]))$ and this last set is compact
since $f$ is proper. Since every finite family of closed sets
$\{\cl\,\phi([t_i,b[\,)\}_i$ has nonempty intersection, then
\[\omega_{\phi}:=\bigcap_{t\in [0,b[}\;
\cl\,\phi([t,b[\,)\ne \emptyset\,,\]
a contradiction.

\hfill \qed

\goodbreak\bigskip

{\bf Closed local homeomorphisms}. The hypothesis of
properness of $f$ can be replaced by closedness of $f$: that is, a
local homeomorphism between Hausdorff spaces which maps closed
subsets of $X$ into closed subsets of $Y$ has the path lifting
property. To see this, argue as above: to prove that $\omega_{\phi}$
is non-empty, take a sequence $t_n\in[0,b[$ converging to $b$ and
such that $\alpha(t_n)$ consists of distinct points,
and is never equal to $\alpha(b)$ (such a sequence certainly
exists, unless $\alpha$ is constant on some left neighborhood
of $b$ ). If the sequence 
$(\phi(t_n))$ has no cluster point, then its range
$R=\{\phi(t_n):n\in\N\}$ is a closed set in $X$; but then
$\{\alpha(t_n):n\in\N\}=f(R)$ is closed in $Y$; this is
plainly absurd, since $\alpha(b)\notin f(R)$, but $(\alpha(t_n))$
converges to $\alpha(b)$. There are relations between
properness and closedness, see Proposition 1.1 below.  
\medskip\goodbreak

{\bf A counterexample}. We are going to show that the
preceding theorem is not true if we drop the Hausdorff
property.
Let $S=\R\cup \{c\}$ with $c\notin \R$ with the following
topology: the open sets in $\R$, $\{c\}\cup A$,
 with $A$ open neighbourhood of $0$
in $\R$, and $\{c\}\cup A\setminus\{0\}$. The topological space
$S$ can be said `the line with two origins',  it is path
connected but the Hausdorff property does not hold true.  We
easily check that the function $f:S\to \R$ whose restriction to
$\R$ is the identity, and with $f(c)=0$, is a proper local
homeomorphism but it is not injective.

Incidentally, also simple connectedness is essential, at least 
for locally well behaved spaces.

\medskip\goodbreak
{\bf Proper maps} Now, let us state two Propositions, whose
proofs are easy, to remind what  
proper functions are in the  context of metrizable spaces
and for maps $\R^n\to \R^m$.

\goodbreak\bigskip

\begin{Proposition}  (Proper maps in metrizable
spaces).   Let $f:X\to Y$ be a continuous function between
the metrizable spaces $X,Y$. Then $f$ is proper if and only if
every sequence $(x_n)$ in $X$ admits a converging
subsequence whenever $(f(x_n))$ converges.
Moreover, if such a function $f$ is proper then it is
closed. 
Finally,  a
closed local homeomorphism between metrizable spaces without
isolated points is a proper map.\end{Proposition}
\goodbreak\bigskip

\begin{Proposition} (Proper maps between
Euclidean spaces). A continuous function $f:\R^n\to \R^m$ is
proper if and only if it is coercive, namely
\[|f(x)|\to \infty\,,\qquad \hbox{{\sl as}}\qquad
|x|\to\infty\,.\]\end{Proposition}

\goodbreak\bigskip
Finally, let us see Gordon's application of the
Hadamard-Caccioppoli Theorem to Algebra. We give some more
details than the original paper [14].

\goodbreak\bigskip

\begin{Proposition} (Nonexistence of a
product in $n$-space for $n\ge 3$).  The $n$-space $\R^n$ with
$n\ge 3$ cannot be endowed of a product operation
$\R^n\times\R^n\to \R^n$, $(x,y)\mapsto xy$ which has the
following properties for any $x,y,z\in \R^n$ and any
$a\in\R$ \begin{itemize}

\item (i)\qquad
$x\,(a\,y)=(a\,x)\,y=a\,x\,y\,,$

\item{} (ii)\quad $\ \;x\,(y+z)=x\,y+x\,z\,,$

\item{} (iii)\quad $\ x\,y=0 \quad \Longrightarrow\quad$ 
either $\ x=0\ $ or $\ y=0\,,$

\item{} (iv)\quad $\
x\,y=y\,x\,$.\end{itemize}\end{Proposition}

\bigskip\goodbreak

\noindent In other words: $\R^n$, with $n\ge 3$, does not
have a commutative  algebra structure without zero
divisors. Remark  that the associative property $\,
x\,(y\,z)=(x\,y)\,z\,$ is not required.

\goodbreak\bigskip

{\it Proof}.
We again argue by contradiction, and we denote by
$F:(x,y)\mapsto xy$ the product. Consider the function $f:X\to
Y$, $x\mapsto x^2=F(x,x)$,  with $X=Y=\R^n\setminus\{0\}$.
First, note that $f$ is a $C^\infty$ function on
$X$: if $x=\sum_{k=1}^nx_ke_k$, where $e_1,\dots,e_n$ is the
standard base of $\R^n$, then
$f(x)=\sum_{k,l=1}^nx_kx_lF(e_k,e_l)$, a quadratic
polynomial function, hence $C^\infty$. Next, denoting by $m,M$
the minimum, respectively the maximum, value of $|f(x)|$ when
$x$ ranges over the unit sphere of $\R^n$, we have
\[0<m|x|^2\leq|f(x)|\leq M|x|^2,\qquad\hbox{{\rm for
every}}\quad x\in X=\R^n\setminus\{0\}\] this follows  from
$|f(x)|=|f(|x|(x/|x|))|=|x|^2|f(x/|x|)|$, valid for every
$x\in
X$ (note that,  by (i), $f(tx)=t^2f(x)$ for every non-zero
real number $t$), and readily implies that $f$ is a proper map.
The differential of $f$ 
 is given by $df(x)v=2xv$, for every $x\in X$ and
$v\in\R^n$. In fact, by (i)  and (ii),
\[f(x+tv)-f(x)=xx+txv+tvx+t^2vv-xx=txv+tvx+t^2f(v)\,;\] 
by (iv)
we then have $f(x+tv)-f(x)=2txv+t^2f(v)$, so that 
\[\lim_{t\to 0}\; {f(x+tv)-f(x)\over t}=2xv+\lim_{t\to
0}\; (tf(v)) =2x\, v.\]
By (iii), $xv=0,x\ne 0$ imply $v=0$. Thus $df(x)$ is nonsingular,
for every $x\in X$. Now all the hypotheses of the
Hadamard-Caccioppoli theorem are satisfied  (in particular $Y$
is simply connected), and so $f$ is a homeomorphism, in
particular it is injective; but clearly $f(x)=f(-x)$, a
contradiction.
 
\hfill \qed

\bigskip

{\it Remark to the proof}.  $Y=\R^n\setminus\{0\}$ is
simply connected if and only if $n\geq3$, and actually commutative 
division
algebra structures exist on $\R^n$
if $n\leq2$; the quaternions prove that commutativity is
essential for the above result (what fails is that $df(x)$, now
given by $df(x)v=xv+vx$, is  singular for some $x\in X$).

\section{Star-shaped images}
In this Section $X,Y$ will always be Banach  spaces,
  $D$ an open connected set, with  $\emptyset\ne D\subseteq X$,
  and $f: D\to Y$ a local homeomorphism.

  \medskip\goodbreak{\bf The auxiliary flow}. Let $x_0\in D$, and
  $y_0=f(x_0)$. 
  We
  are going to define a flow 
  $\Phi: \DomPhi\to D$ which will be our tool in investigating
  the invertibility of $f$ around $x_0$. The basic properties
  of $\Phi$, so that it is called a flow in~$D$, are the following:
  \begin{itemize}
  \item (i) $D_\Phi$ is an open subset of~$D\times\R$, and
  $\Phi\colon D_\Phi\to D$ is continuous; 
  
  \item (ii) for all $x\in D$, the set
  $\{t\in\R\;:\;(x,t)\in D_\Phi\}$ is an interval
  containing~$0$;
  
  \item (iii) $\Phi(x,0)=x\quad\hbox{for all }x\in D\,$;
  
  \item (iv) $(x,t_1),(x,t_1+t_2)\in
  D_\Phi\;\Rightarrow\;
  \bigl(\Phi(x,t_1),t_2\bigr)\in D_\Phi\;$  and 
  $\Phi\bigl(\Phi(x,t_1),t_2\bigr)=\Phi(x,t_1+t_2)\;$
  for all $x\in D$, $t_1,t_2\in\R\,$.\end{itemize}
 
  \noindent If $\{x\}\times[0,+\infty[\subset \DomPhi$ we say
  that the trajectory through $x$ is {\em global in the future}.
  Moreover, whenever  $\DomPhi=D\times\R$ we say
  that $\Phi$ is a (global) {\em dynamical system} in~$D$.
  
  To define $\Phi$ we start from  the following
  dynamical system in~$Y$:
  \begin{equation}\Psi\colon Y\times\R\to Y\,,\qquad
  \Psi(y,t):=y_0+e^{-t}(y-y_0)\,,
  \end{equation}
  whose trajectories are the half-lines hinged at~$y_0$, but 
  with an exponential parameter instead of a linear one, so that
  $\Psi(y,0)=y$, $\Psi(y,t)\to y_0$ as~$t\to+\infty$.
  It is indeed a
  dynamical system, because $\Psi(\Psi(y,t_1),t_2)=\Psi(y,t_1+t_2)$.
  \goodbreak\bigskip
  
  \begin{Lemma} (The auxiliary flow). 
  Let $X,Y$ be Banach spaces, let $D\subseteq X$ be 
  open and connected, let $x_0\in D$, and let $f:D\to Y$
  be a local homeomorphism. Then there exists a flow 
  $\Phi:\DomPhi\to D$
  which satisfies the following formula

  \begin{equation}f\bigl(\Phi(x,t)\bigr)=\Psi\bigl(f(x),t\bigr)\quad
  \hbox{for all }(x,t)\in D_\Phi\,,\end{equation}
  and two such flows coincide in the intersection of their domains
  (so  $\Phi$ will be maximal in the sequel). In the particular
  case where $f$ is a local diffeomorphism
  (namely it is also $C^1$ together
  with all its local inverses),  the mapping $\Phi$ is $C^1$ and
  it is the flow of the following differential equation
  \begin{equation}\dot x=F(x)\,,\qquad F: D\to X\,,\  x\mapsto
  -f'(x)^{-1}\left(f(x)-f(x_0)\right)\,.
\end{equation}\end{Lemma}

  \bigskip
  
  \noindent In other words we could say that $\Phi$ is the
  maximal
  lifting
  of $\Psi\circ (f\times id)$ (where $id$ the
  identity in $\R$) such that $\Phi(x,0)=x$ for all $x\in D$.
  
  \bigskip

  {\it Proof.} Fix $x\in D$ and consider the continuous function
  $\;\R\to Y\,,$ $\;t\mapsto \Psi(f(x),t)\,$.
  By similar arguments as in Section 1 we prove the existence of a
  unique maximal lifting $]a(x),b(x)[\to D$, $t\mapsto
  \Phi(x,t)$,
  with $\Phi(x,0)=x$, $-\infty \le a(x)<0<b(x)\le +\infty$.
  Let $\DomPhi:=\bigcup_{x\in D}\{x\}\times ]a(x),b(x)[$. All
  the properties above are easy to check except (i) which
  requires some arguments. 
   
   We consider the subset $D\times
  [0,+\infty[$ only; the set $D\times ]-\infty,0]$ is
  handled similarly.
  Let $x_0\in D$ be given. 
  First consider the supremum $\tau (x_0)$ of all
  real numbers $t\geq0$ such that
  $\{x_0\}\times[0,t[$ is contained in the
  interior of $\DomPhi$ (if no such $t>0$ exists,
  then $\tau (x_0)=0$). Next, define $E$ to be the set
  of all real  $t\in[0,\tau (x_0)]$ such that $\Phi$ is
  not continuous at $(x_0,t)$; arguing as in Lemma
  1.2 one easily sees that $E$ is empty.
  And still arguing as in Lemma 1.2, with $\tau (x_0)$
   in place of $a$, it is also easy to see that
   $\tau (x_0)=b(x_0)$, hence that $\DomPhi$ is open.

   \hfill\qed 

   \goodbreak\bigskip

   {\bf The attraction basin}. Let $f:D\to Y$, $x_0$,
   $\Phi$ be as in Lemma 2.1 (in the general case), and let 
   $y_0=f(x_0)$.       
   Let $U$ be an open neighbourhood of~$x_0$ where $f$ is
   injective and let $g:=(f|U)^{-1}$. For any small $r>0$,
   the ball $B(y_0;r)$ (with center at $y_0$ and radius $r$) 
   is contained in~$f(U)$, and for
   such~$r$ let $U_r:= g^{-1}(B(y_0;r))$. Then $U_r$ is a
   neighbourhood of~$x_0$, and for all~$x\in U_r$ the
   trajectories $t\mapsto\Phi(x,t)$ of~$\Phi$ are defined
   globally in the future, belong to~$U_r$ for all~$t\ge0$
   and converge to~$x_0$ as~$t\to+\infty$. Then $x_0$ is
   an {\it attractor} namely it attracts a whole neighbourhood
   (any $U_r$ will do), and it is {\it stable} that is any of its
   neighbourhoods contains a positively invariant
   neighbourhood with
   global existence in the future, indeed
   again we can consider $U_r$, with small enough~$r$
   (we remind that positive invariance means that
   $\Phi(x,t)\in U_r$ for any $x\in U_r$ and $t>0$, such that
   $(x,t)\in\DomPhi$). So
   we just
   proved that
   $x_0$ is  asymptotically stable, i.e. a stable attractor.
   
   The maximal neighbourhood  $\A$ of $x_0$ such that,
   for all~$x\in \A$, the
   trajectories $t\mapsto\Phi(x,t)$ of~$\Phi$ are defined
   globally in the future, belong to~$\A$ for all~$t\ge0$,
   and converge to~$x_0$ as~$t\to+\infty$, is called the
   basin of attraction of~$x_0$.
   
   \bigskip

   \begin{Proposition}(Injectivity in the
   attraction basin). Under the hypotheses of the first part of
   Lem\-ma~2.1 the
   attraction basin $\A$ of $x_0$ for $\Phi$ is open. Moreover:
   \begin{itemize}
   \item (i) the
   restriction of $f$ to  $\A$ is
   injective, 
   \item (ii) $f(\A)$ is star-shaped with respect to
   $y_0:=f(x_0)$, and 
   \item (iii) $\A$ is the maximal connected subset of
   $D$ which contains $x_0$ and has the properties (i) and
   (ii).\end{itemize}\end{Proposition}
   
   \bigskip\goodbreak

   {\it Proof.} $\A$ is open because $\DomPhi$ is open in
   $X\times\R$ and $\Phi$ is continuous. 
   
   To prove that $f$ is
   injective on $\A$, let $x_1,x_2\in \A$ be such that
   $f(x_1)=f(x_2)$.
   Then for all $t\ge 0$ 
   \[f(\Phi(x_1,t))=y_0+e^{-t}(f(x_1)-y_0)=
   y_0+e^{-t}(f(x_2)-y_0)=f(\Phi(x_2,t))\,.\]
   Since, for large $t$, both $\Phi(x_1,t)$ and $\Phi(x_2,t)$
   enter a neighbourhood of $x_0$ where $f$ is injective, we
   have that  $\; \Phi(x_1,t)=\Phi(x_2,t)$  for large
   $t$. Thus for large $t$ we have 
   $\; x_1=\Phi(\Phi(x_1,t),-t)=\Phi(\Phi(x_2,t),-t)=x_2\,$.
   The image $f(\A)$ is star-shaped with respect to
   $y_0$ because
   $$f(\A)=\{y_0\}\cup \{y_0+e^{-t}(f(x)-y_0):
   (x,t)\in\DomPhi\}\,.$$

   The maximality is also easily verified.
   
   \hfill\qed 

\goodbreak\bigskip
  
    \begin{Proposition}(Bijectivity $\;
\Longleftrightarrow\;$ 
$\DomPhi=D\times\R$). 
Let $X$ and $Y$  be 
Banach spaces, let $D\subseteq X$ be open and connected, let
$x_0\in D$, let $f:D\to Y$ be a local homeomorphism,
and let $\Phi$
be the auxiliary flow as above. Then $f$ is a global
homeomorphism onto $Y$ if and only if the flow $\Phi$ is a
global dynamical system.\end{Proposition}\goodbreak\bigskip

    {\it Proof.} Suppose first that $f$~is a global
homeomorphism onto~$Y$. Then the inverse mapping $f^{-1}$~is
defined and continuous on~$Y$ and the expression $\Phi(x,t)=
f^{-1}(y_0+e^{-t}(f(x)-y_0))$ is defined and continuous for
all~$(x,t)\in D\times\R$.

\noindent
Conversely, suppose that $\DomPhi=D\times\R$. Let~$y\in\Rn$
and~$\varepsilon>0$ such that $y_0+\varepsilon(y-y_0)\in
f(\A)$, and let $g:=(f|\A)^{-1}$. Then
\[ f(D)\supset f(\A)\ni f\Bigl(\Phi\bigl(g(y_0+
  \varepsilon(y-y_0)),\ln\varepsilon\bigl)\Bigl)=
  y_0+e^{-\ln\varepsilon}
  \bigl(y_0+\varepsilon(y-y_0)-y_0\bigr)=
  y\,,
  \]
and $f|\A$~is proved to be onto~$Y$. To verify that $f$~is
also one-to-one on all of~$D$, i.e., that $\A=D$, it
suffices to prove that $\A$ is a closed subset of~$D$,
because we already know that it is open and nonempty. Let
then $x_n\in\A$ be a sequence converging to~$x\in D$. Since
$f(\A)=Y$, there exists $\bar x\in\A$ such that $f(\bar
x)=f(x)$. Recalling that $(f|\A)^{-1}: Y\to\A$ is continuous,
from $f(x_n)\to f(\bar x)$ we get that $x_n\to\bar x$,
whence $x=\bar x\in\A$.

\hfill\qed
    
\goodbreak\bigskip
Now, let us prove Theorem~0.2 in the Introduction.
\goodbreak\bigskip

{\bf Proof of the Hadamard-Levy Theorem}. By the preceding
Proposition 2.3 we can just show that the solutions to the
equation (2.3) are defined on the whole $\R$. First remark
that by (2.1), and (2.2),
\[\|f(\Phi(\bar x,t))-y_0\|=e^{-t}\|f(\bar x)-y_0\|\]
so this is bounded whenever $t$ ranges on a bounded interval.
Then, along a trajectory
$\gamma:]a,b[\to D$, $\gamma(t)=\Phi(\bar x,t)$,  defined in
a bounded interval of time $]a,b[$, we have the following
estimate for the vector field in (2.3):
$$\|F(\gamma(t))\|\le
\|f'(\gamma(t))^{-1}\|\; \|f(\gamma(t))-y_0\|\le
c\;\beta(\|\gamma(t)\|)\,,$$ 
for a suitable $c>0$ (the function
$\beta$ was introduced in (0.3)). 

From now on the arguments are standard, however we prefer
to complete the proof to be self-contained. Let $r(t):=
\|\gamma(t)\|$. Then
for $a\le t_1\le t_2\le b$ we have
\begin{equation}\|r(t_2)-r(t_1)\|\le
  \|\gamma(t_2)-\gamma(t_1)\|\le c\int_{t_1}^{t_2}\;
  \beta\left(\|\gamma(t)\|\right)dt\,.
  \end{equation}
The function $x\mapsto \|x\|$ is Lipschitz continuous and the
function $\gamma$ is~$C^1$ (remind that $f$ is a local
diffeomorphism in the present theorem), so that $t\mapsto
r(t)$ is locally absolutely continuous and it has derivative
almost everywhere. By the previous estimate, dividing by
$t_2-t_1$ and going to the limit we have $\|r'(t)\|\le c\;
\beta (r(t))$ almost everywhere. Now, for $t,t_0\in ]a,b[$
\[\left|\int_{r(t_0)}^{r(t)}\; {1\over{\beta(s)}}\,
  ds\right|=\left|\int_{t_0}^t\;{r'(s)\over{\beta(r(s))}}\,
  ds\right|\le 
  \left|\int_{t_0}^t\;\Bigl|{r'(s)\over{\beta(r(s))}}
  \Bigr|\,ds\right|\le c\;|t-t_0|
  \le c\;|b-a|\,.\]
Then $r(t)$ for $t\in]a,b[$ is bounded from above by any
$r_0>0$ large enough to give 
$\int_{r(t_0)}^{r_0}\;{1\over\beta(s)}\;ds \ge c\;|b-a|$
(remind the first formula in (0.3)).  Using again the
inequality~(2.4) and this time the monotonicity of~$\beta$ we
see that $\|\gamma'(t)\|\le c\;\beta(r_0)$. Then $\gamma$~is
Lipschitz continuous on~$]a,b[$ and it can be extended by
continuity to~$a$ and~$b$.

\hfill\qed

   \goodbreak\bigskip
   
   In the sequel we shall need the following Lemma:
   \goodbreak\bigskip

   \begin{Lemma} (On $\partial \A$ the
   trajectories have finite life). Let us assume the hypotheses of
   the first part of Lemma 2.1. Then the attraction basin $\A$
   is invariant, namely
   $\;x\in\A\ $ $ \Longrightarrow$ $ \ \Phi(x,t)\in \A\;$ 
   for all
   $t$ such that $(x,t)\in\DomPhi$, and  also 
   $\partial\A$ (the boundary of $\A$ in $D$) is invariant.
   Moreover,
   there is not global existence in the future
   for~$t\mapsto \Phi(x,t)$ if
   $x\in\partial\A$.\end{Lemma}

   \bigbreak

   {\it Proof.} First of all let us see that
  $$f(\partial\A)\subseteq\partial f(\A)\,.
  \eqno(2.5) $$
  The set $\A$ is open in~$X$ and $f$ is a one-to-one local
  homeomorphism on~$\A$, so that $f(\A)$ turns out to be
  open, too, and $f\vert \A\colon\A\to f(\A)$ is a
  homeomorphism. $f(\partial\A)$ is contained
  in the closure of~$f(\A)$ because $f$ is continuous.
  Let $\bar x$ be a point in the closure of~$\A$ such
  that $f(\bar x)\in f(\A)$, i.e., $f(\bar x)=
  f(x)$ for some~$x\in\A$. Let~$x_n$, $n\ge1$, be a
  sequence of points of~$\A$ converging to~$\bar x$.
  By continuity of~$f$ we have
  $f(x_n)\to f(\bar x)=f(x)$,
  and by continuity of~$(f\vert \A)^{-1}$ we have
  $x_n=(f\vert \A)^{-1}(f(x_n))\;\to\;
  (f\vert \A)^{-1}(f(x))=x\,$,
  so that $\bar x=x\in\A$.
  From~(2.5) and the fact that $f(\A)$ is a neighbourhood
  of~$y_0=f(x_0)$, there exists $\varepsilon>0$ such that
  
  \begin{equation}
  x\in\partial\A\qquad\Longrightarrow\qquad
   \|f(x)-y_0\|\ge\varepsilon
  \end{equation}
  
  It is obvious from its definition that $\A$ is invariant
  for the flow $(x,t)\mapsto \Phi(x,t)$. 
  The same holds for $\partial\A$:
  In fact, let $x\in\partial\A$, $x_n\in\A$,
  $x_n\to x$, $(x,t)\in \DomPhi$. Then $(x_n,t)\in \DomPhi$ for
  all
  large~$n$, because $\DomPhi$ is open, and, by continuity
  $\A\ni \Phi(x_n,t)\to \Phi(x,t)$.  
  The point $\Phi(x,t)$ belongs to the closure
  of~$\A$, but not to~$\A$, because otherwise
  $x= \Phi(\Phi(x,t),-t)$ itself would be in~$\A$.

  Finally, from~(2.6) we get:
  \[x\in\partial\A\,\,
  \Longrightarrow\,\,
  \varepsilon\le \|f(\Phi(x,t))-y_0\|=
  e^{-t}\|f(x)-y_0\|\,\,\Longrightarrow\,\,
  t\le\ln{\|f(x)-y_0\|\over\varepsilon}\,.
  \]

  \hfill \qed
  
   \bigskip

    \goodbreak{\bf Bounded sets in $D$}. In the sequel we say that
    a set $B\subseteq D$ is bounded in $D$ if (i) it is
    bounded as a subset of $X$, and (ii) its
    closure in $X$ is contained in $D$. 
    
    \goodbreak\medskip
    
    {\bf Trapped trajectories}. We
    need to guarantee that the trajectories
    of $\Phi$ which are trapped into a closed and bounded subset
    of $D$ are defined globally in the future (condition (c) in
    Lemma 2.3
    below).  This  is familiar and always true
    for solutions to differential equations  which  are  `trapped' into
    compact sets in
    finite dimension. The following Lemma 2.3 shows few
    technical conditions each of which
    implies this property. In the  statement
    we denote by $\;[f(x_0);f(x)]\subset Y\;$  the line
    segment from $f(x_0)$ to $f(x)$.
    
    \bigskip

    \begin{Lemma}(Trapped
    trajectories never die). 
    Let $X$ and $Y$  be 
Banach spaces, let $D\subseteq X$ be open and connected, let
$x_0\in D$, let $f:D\to Y$ be a local homeomorphism,
and let $\Phi$
be the auxiliary flow as above. Consider the following
conditions:
\begin{itemize}

 \item (a-1) the restriction $f|B$ is proper for any 
set $B$ closed and bounded in $D$;

 \item (a-2) $f$ is a local $C^1$ diffeomorphism and
for each bounded and closed set $B\subset D$ we have
\begin{equation}\sup_{x\in B}\|f'(x)^{-1}\|<+\infty\end{equation}

 \item (b)  for any $B$, closed and bounded subset of $D$, and any
$x\in B$, the connected components of
$f^\leftarrow([f(x_0);f(x)])\cap B$ are compact; 

 \item (c) for any
$B$, closed and bounded subset of $D$, and any $x\in B$, if
$\Phi(x,t)\in B$ for all $t>0$ such that $(x,t)\in\DomPhi$, then
the trajectory through $x$ is global in the future (in other words:
trajectories which are eventually in bounded closed sets never
die).\end{itemize}
 \medskip
 
 Then either one of (a-1) and (a-2) imply (b),
which implies (c). All conditions are trivially
satisfied if $X$ is finite dimensional.\end{Lemma}
\goodbreak\bigskip

    {\it Proof.} 
   The proof is trivial except for 
 (a-2) $\Rightarrow$ (b). Let $L$ be a component
of $f^\leftarrow([f(x_0);f(x)])\cap B$. Pick $\, x_1\in L$, and
let $v=f(x)-f(x_0)$. If $v=0$, then $[f(x_0);f(x)]$ consists of
the single point $f(x_0)$ and $L$ is then also a singleton, since
$f$ is a local homeomorphism. Assume then $v\not=0$. Since $f(L)$
is connected, the set $\{t\in\R:f(x_1)+tv\in f(L)\}$ is a
bounded interval $I$ of $\R$ containing $0$. Let
$\alpha:J\to L$ be the maximal lifting of the path 
$\ell(t)=f(x_1)+tv\,(t\in I)$ with origin
$\alpha(0)=x_1$. Since $f$ is a local diffeomorphism,
such an $\alpha$ is differentiable, and
differentiating $f(\alpha(t))=f(x_1)+tv$ we get
$f'(\alpha(t))(\alpha'(t))=v$, whence
$\alpha'(t)=f'(\alpha(t))^{-1}v$. Since $\sup_{x\in
L}\|f'(x)\|$ is finite, $\alpha'$ is bounded on its maximal
interval $J$ of existence; thus the $\omega$-limit set of $\alpha$
is nonempty, and it is contained in the closed set
$L$. It follows that $J=I$, and by the same token, that $\inf
I\in I$, and $\sup I\in I$, that is, $I$ is compact. It is now
obvious that $f$ induces a homeomorphism of $L$ onto $f(L)$,
which has $\alpha\circ\ell^{-1}$ as inverse. Thus $L$ is compact,
since $f(L)$ is homeomorphic to $I$ {\sl via} $\ell$.

\hfill\qed
    
    \goodbreak\bigskip
    
    {\bf A class of functions satisfying (a-1)}.
    The condition (a-1) is fulfilled  
    if $f=p+c$ with $p$  proper  and $c$  compact,
    i.e., mapping closed bounded sets to compact sets. Indeed,
    remind Proposition 1.1, and consider a sequence
    $(x_n)$ in the closed bounded set $B$, with $(f(x_n))$
    convergent. Since $c$ is
    compact, it maps  a
    subsequence $(x_{n_k})$ to a convergent sequence
    $(c(x_{n_k}))$, thus $p(x_{n_k})=f(x_{n_k})-c(x_{n_k})$
    converges
    and finally $(x_{n_k})$ has a convergent subsequence since
    $p$ is proper.
    
    \goodbreak\bigskip

    {\bf Coercive auxiliary functions}. 
    The nonnegative continuous function
    $k:D\to \R$ is called  coercive whenever for
    any $a>0$  the inverse
    image $k^{\leftarrow}([0,a])$ is
    bounded in $D$.
     
    \goodbreak\medskip
     
    {\bf Global Lyapunov functions}.   In our
    framework  the 
    function   $k:D\to \R_+$ is called a global Lyapunov function
    for the flow $\Phi$  above, if it is continuous,
    nonnegative, coercive, and
    weakly decreasing along the trajectories,
    namely $t\mapsto
    k(\Phi(x,t))$ weakly decreases for all $x\in D$.

  \begin{Proposition} (Injectivity and 
    star-shaped image by
    Lyapunov functions). 
    Let $X,Y$ be Banach spaces, let $D\subseteq X$ be 
    open and connected, let $x_0\in D$, and let $f:D\to Y$
    be a local homeomorphism. 
    Then  $f$ is injective, and  the image $f(D)$ is
    star-shaped with respect to
    $f(x_0)$, if  there exists a global Lyapunov
    function for $\Phi$, and $f$ satisfies any of the
    conditions (a-1), (a-2), (b), (c) in Lemma~2.3. 
    \end{Proposition}
   
    \bigskip\goodbreak

    {\it Proof.}  We are
    going to prove that $D=\A$.  So  we are done by
    Proposition~2.1.
     
    It is enough to show that the boundary
    $\partial \A$ (of $\A$ in $D$) is empty.    
    We argue by contradiction and assume that $x\in \partial\A$
    By Lemma~2.2 the maximal positive trajectory through $x$,
    $\gamma:[0,b[\to D$, $t\mapsto \Phi(x,t)$, lies in
    $\partial\A$, and has a finite life:
    $\;\gamma([0,b[\,)\subseteq \partial\A$, and $b<+\infty$.

    The Lyapunov function  $k:D\to\R_+$
    is coercive and, in particular,
    $B:=k^{\leftarrow}([0,b])$ is bounded in $D$. Moreover,
    $k\circ \gamma$ is decreasing and so
    $\;\gamma([0,b[\,)\subseteq B$. Now Lemma 2.3 says that
    condition (c) above holds true, namely $b=+\infty$, a
    contradiction.
    
    \hfill\qed

      \goodbreak\bigskip
 
    The preceding result, as well as the following one, 
extend some  results in [17] (by two of the authors) 
 where the  finite dimensional case is treated. 
    That paper also shows that the converse of Proposition~2.3
    holds true  in $\R^n$ (and proves other related facts).
In the following statement we consider an Hilbert space $X$
with scalar product `$\,\cdot \,$', and $B(x_0;r)$ will denote
the open ball $\|x-x_0\|<r$. We could formulate an analogous
 fact in general Banach spaces but it would be more
complicated to be stated (but not to be proved).

     \goodbreak\bigskip

    \begin{Proposition}(A criterion of
injectivity on a ball). 
Let $X$ be a Hilbert space, $x_0\in X$, $Y$ be a Banach space,
$f\colon B(x_0;r_0)\to Y$ be a local $C^1$~diffeomorphism
satisfying any of the conditions of Lemma~2.3. Then the
following two conditions are equivalent:\begin{itemize}

\item (a) $f$~is injective and $f(B(x_0;r))$ is
star-shaped with respect to~$f(x_0)$ for all positive
 $r\le r_0$;

\item (b) the following inequality holds 
for all $x\in B(x_0;r_0)$
\begin{equation}(x-x_0)\cdot f'(x)^{-1}\bigl(f(x)-f(x_0)\bigr)\ge0
  \end{equation}

\end{itemize}\end{Proposition}
\bigskip\goodbreak

{\it Proof.} The left-hand side of~(2.8) is the derivative
with respect to~$t$ at~$t=0$ of the scalar function
$$t\mapsto{1\over2}\|\Phi(x,t)-x_0\|^2\,.
 $$
Asking it to be nonnegative is the same as asking the scalar
function $x\mapsto(1/2)\|x-x_0\|^2$ to be weakly decreasing
along the flow~$\Phi$, which in turn is the same
as requiring the same from each of the functions $x\mapsto
1/(r^2-\|x-x_0\|^2)$ on~$B(x_0;r)$, $0<r\le r_0$. These last
functions have the advantage of being coercive
on~$B(x_0;r)$. Hence condition~(b) is satisfied,
Proposition~2.6 can be applied to get condition~(a).

\noindent
Conversely, if condition~(a) holds, then the sets $B(x_0;r)$
are positively invariant for~$\Phi$ and the (square) norm
of~$\Phi(x,t)$ must be a weakly decreasing function of~$t$,
whence inequality~(2.8).

\hfill\qed


\begin{thebibliography}{99}

\bibitem{1} A. AMBROSETTI and G. PRODI,
On the inversion of some differentiable mappings
with singularities
    between Banach spaces,
{\it Ann. Mat. Pura Appl.} {\bf 93}, 231--247 (1973).

\smallbreak

\bibitem{2}
 S. BANACH and S. MAZUR,
 \"Uber mehrdeutige stetige Abbildungen,
{\it Studia Math.}{\bf~5}, 174--178 (1934).

\smallbreak

\bibitem{3}
 M.S. BERGER,
 {\it Nonlinearity and functional analysis},
  Academic Press, 1977.

\smallbreak

\bibitem{4}
 N.P. BHATIA and G.P. SZEG\"O,
 {\it Stability theory of dynamical systems},
 Springer-Verlag, 1970.

\smallbreak

\bibitem{5}
 F. BROWDER,
 Covering spaces, fiber spaces and local homeomorphisms,
{\it Duke Math. J.}{\bf~21}, 329--336 (1954).

\smallbreak

\bibitem{6}
K.J. BROWN and S.S. LIN,
 Periodically perturbed conservative systems and a global
inverse function theorem,
{\it Nonlinear Analysis TMA} {\bf 4} , 193--201 (1980).

\smallbreak

\bibitem{7}
 R. CACCIOPPOLI,
 Sugli elementi uniti delle trasformazioni funzionali,
{\it Rend. Sem. Mat. Univ. Padova}~{\bf3}
    , 1--15 (1932).

\smallbreak

\bibitem{8}
 R. CACCIOPPOLI,  Un principio di inversione per 
le corrispondenze funzionali e sue applicazioni alle
equazioni alle derivate parziali, {\it Atti Acc. Naz. 
Lincei}~{\bf16}, 390--400 (1932).

\smallbreak

\bibitem{9}
 S.N. CHOW and  J.K. HALE, {\it  Methods of bifurcation
theory}, Springer-Verlag, 1982.

\smallbreak

\bibitem{10}
 L. M. DRU\.ZKOWSKI and H.K. TUTAI,
 Differential conditions to verify the Jacobian
conjecture\rm, 
{\it Ann. Polon. Math.} {\bf 57}, 253--263 (1992).

\smallbreak

\bibitem{11}
\rm D. GALE and H. NIKAIDO,
 The Jacobian matrix and global univalence of mappings,
{\it Math. Ann.}~{\bf 159}, 81--93 (1965).

\smallbreak

\bibitem{12} W.B. GORDON,
 On the diffeomorphisms of Euclidean space\rm,
{\it Amer. Math. Monthly}~{\bf 79}, 755--759 (1972).

\smallbreak

\bibitem{13}  W.B. GORDON,
 Addendum  to ``On the diffeomorphisms of Euclidean
       space'',
 {\it Amer. Math. Monthly}~{\bf 80}, 674--675 (1973).

\smallbreak

\bibitem{14}
W.B. GORDON, An application of Hadamard's inverse function
theorem to algebra,
{\it Amer. Math. Monthly}\ {\bf 84}, 28--29 (1977).

\smallbreak

\item{15}
\rm G. GORNI,
 A criterion of invertibility in the large for local 
diffeomorphisms
    between Banach spaces,
{\it  Nonlinear Analysis~TMA} {\bf 21}, (1993) 43--47.
    

\smallbreak

\bibitem{16} G. GORNI and G. ZAMPIERI,   
Global sinks for planar vector fields, 
{\it Evolution Equations and Nonlinear Problems, 
Proceedings
of the RIMS Symposium, RIMS Ko\-kyu\-ro\-ku} {\bf 785}, 
Ky\=oto,
134--138 (1992).

\bibitem{17}    G. GORNI, and G. ZAMPIERI,  
 Injectivity onto a star-shaped set
for local homeomorphisms in n-space,
{\it Annales Polonici
Mathematici} {\bf 59},  171--196  (1994).

\smallbreak

\bibitem{18} C. GUTIERREZ,
 Dissipative vector fields on the plane with infinitely
many attracting hyperbolic singularities,
{\it Bol. Soc. Bras. Mat.} {\bf 22}, 179--190 (1992).

\smallbreak

\bibitem{19}
 J. HADAMARD,
 Sur les transformations ponctuelles,
{\it Bull. Soc. Math. France} {\bf 34}, 71--84 (1906).

\smallbreak

 \bibitem{20} J. HADAMARD,  Sur
 les correspondances ponctuelles. 
 {\it Oeuvres~I}, Editions du CNRS,
 383--384 (1968).
 
\smallbreak

\bibitem{21} S. INVERNIZZI and F. ZANOLIN, On
the existence and uniqueness of periodic solutions
of differential delay equations,
{\it Math. Z.} {\bf 163}, 25--37 (1978).

\smallbreak



\bibitem{22}
 M. P. LEVY,
Sur le fonctions de ligne implicites,
{\it Bull. Soc. Math. France} {\bf 48}, 13--27 (1920).

\smallbreak

\bibitem{23}
 G.H. MEISTERS,
 Inverting polynomial maps of $n$-space by solving 
differential equations\rm,
in {\it  Fink, Miller, Kliemann Editors,  Delay and Differential
Equations: Proceedings in Honour of George Seifert on his
retirement}, World Sci. Pub.  Co., 107--166 (1992).

\smallbreak

\bibitem{24}
 G.H. MEISTERS and C. OLECH,
 Locally one-to-one mappings and a classical 
theorem on schlicht
    functions,
   {\it Duke Math. J.} {\bf 30}, 63--80 (1963).

\smallbreak

\bibitem{25}
 G.H. MEISTERS and C. OLECH,
 Solution of the global asymptotic stability 
Jacobian conjecture for the
    polynomial case,
in: {\it Analyse Math\'ematique et applications}, 
Gauthier-Villars, Paris ,
    373--381 (1988).
 
 \smallbreak
 
\bibitem{26}
 G.H. MEISTERS and C. OLECH,
 Global stability, injectivity, and the Jacobian
conjecture, {\it Proceedings of the first World Congress of
Nonlinear Analysts}, to appear.

\smallbreak

\bibitem{27}
\rm C. OLECH,
On the global stability of
 an autonomous system on the plane,
   {\it Cont. Diff. Eq.} {\bf 1},
    389--400 (1963).

\smallbreak

\bibitem{28}
 J.M. ORTEGA and  W.C. RHEIBOLDT,
{\it Iterative solutions of nonlinear equations
in several variables},
Academic Press, 1970.

\smallbreak

\bibitem{[29]}
\rm T. PARTHASARATHY,
 On global univalence theorems\rm,
   {\it Lecture Notes in Math.\penalty5000\ 977},
   Spring\-er Verlag, 1983.

\smallbreak

\bibitem{30}
\rm R. PLASTOCK,
 Homeomorphisms between Banach spaces\rm,
{\it Trans. Amer. Math. Soc.} {\bf 200}, 169--183 (1974).

\smallbreak

 \bibitem{[31]} G. PRODI and A. AMBROSETTI,
 {\it
 Analisi non lineare}, Quaderni della Scuola
 Normale Superiore, Pisa, Italy
 1973.
 
 \smallbreak

\bibitem{[32]}
\rm P.J. RABIER,
On global diffeomorphisms of Euclidian space,
{\it Technical Report ICMA-91-159},
    Pittsburgh (1991).

\smallbreak

\bibitem{33}
\rm M. RADULESCU and S. RADULESCU,
 Global inversion theorems and applications 
to differential equations\rm,
{\it Nonlinear Analysis TMA} {\bf 4}, 951--965 (1980).

\smallbreak

\bibitem{34}
\rm W. C. RHEINBOLDT,
 Local mapping relations and global
implicit function theorems,
 {\it Trans. Amer. Math. Soc.} {\bf 138}, 183--198 (1969).

\smallbreak

\bibitem{[35]} M. SABATINI, An extension to Hadamard global
inverse function theorem in the plane,
{\it Nonlinear Analysis TMA}, to appear.

\smallbreak

\bibitem{36} I. W. SANDBERG, 
 Global inverse
function theorems, {\it I.E.E.E. Trans. Circuits
Systems CAS} {\bf 27}, 998--1004 (1980).

\smallbreak

\bibitem{37} S. SOLIMINI and C. MARICONDA,
 Note sui teoremi sulla funzione implicita e costruzione
del grado topologico, {\it S.I.S.S.A., Trieste, Italy} (1988).

\smallbreak

\bibitem{38}
\rm J. SOTOMAYOR,
 Inversion of smooth mappings,
{\it Z. Angew. Math. Phys.} {\bf 41}, 306--310 (1990).

\smallbreak

\bibitem{37}
G. VIDOSSICH, Two remarks on the stability of ordinary
differential equations,
{\it Nonlinear Analysis TMA} {\bf  4}, 967--974 (1980).

\smallbreak

\bibitem{40}
T. WA\.ZEWSKI,
 Sur l'evaluation du domain d'existence de 
fonctions implicites
    r\'eelles ou complexes\rm,
{\it Ann. Soc. Polon. Math.} {\bf 20}, 81--120 (1947).

\smallbreak

\bibitem{41}
\rm G. ZAMPIERI,
 Finding domains of invertibility for smooth functions
by means of
    attraction basins\rm,
  {\it  J. Differential Equations} {\bf 104}, 11--19 (1993).

\smallbreak

\bibitem{42}
\rm G. ZAMPIERI,
 Diffeomorphisms with Banach space domains\rm,
{\it    Nonlinear Analysis TMA} {\bf 19},  923--932 (1992).

\smallbreak

\bibitem{43}
G. ZAMPIERI and G. GORNI,
 On the Jacobian conjecture for global asymptotic
stability,
{\it J.  Dynamics Diff. Eq.} {\bf 4}, 43--55 (1992).
 
 \smallbreak
 
 \bibitem{44}  G. ZAMPIERI and G. GORNI,  
 Local homeo- and diffeomorphisms: invertibility and
convex image\rm, {\it Bulletin of the Australian
Mathematical Society} {\bf 49},   377--398  (1994).
\end{thebibliography}
\end{document}